\documentclass[12pt,a4paper]{article}
\usepackage{amssymb, amscd, amsthm, amsmath,latexsym, amstext,mathrsfs}
\usepackage[all]{xy}
\usepackage[dvips]{graphicx}

\newtheorem{thm}{Theorem}[section]
\newtheorem{coro}[thm]{Corollary}
\newtheorem{lemma}[thm]{Lemma}
\newtheorem{prop}[thm]{Proposition}
\newtheorem{conj}[thm]{Conjecture}

\theoremstyle{remark}
\newtheorem{remark}[thm]{\textbf{Remark}}

\theoremstyle{definition}

\newtheorem{example}[thm]{Example}

\numberwithin{equation}{thm}
\newcommand{\set}[1]{\{\,{#1}\,\}}

\newcommand{\oi}{\mathrm{i}}
\newcommand{\ii}{\mathrm{ii}}
\newcommand{\iii}{\mathrm{iii}}

\newcommand{\dgf}[1]{\langle\, #1\,\rangle}

\begin{document}
\title{The Pythagoras number and the $u$-invariant of Laurent series fields in several variables}
\author{Yong HU}
\date{}

\maketitle


\begin{abstract}
We show that every sum of squares in the three-variable Laurent series field $\mathbb{R}(\!(x,y,z)\!)$ is a sum of $4$ squares, as was conjectured in a paper of Choi, Dai, Lam and Reznick in the 1980's. We obtain this result by proving that every sum of squares in a finite extension of $\mathbb{R}(\!(x,y)\!)$ is a sum of $3$ squares. It was already shown in Choi, Dai, Lam and Reznick's paper that every sum of squares in $\mathbb{R}(\!(x,y)\!)$ itself is a sum of two squares. We give a generalization of this result where $\mathbb{R}$ is replaced by an arbitrary real field. Our methods yield similar results about the $u$-invariant of fields of the same type.
\end{abstract}

\section{Introduction}

Let $K$ be a field, which we assume to have characteristic different from 2.
 The Pythagoras number $p(K)$ of  $K$ is the smallest integer $p\ge 1$ or $+\infty$ such that every sum of (finitely many) squares in $K$ can be written as a sum of at most $p$ squares in $K$. (Of course, this definition is still valid in characteristic 2.) The $u$-invariant $u(K)$  in the sense of Elman--Lam \cite{ElmanLam73MZ} is the supremum of dimensions of anisotropic \emph{torsion} quadratic forms over $K$. (A quadratic form over $K$ is called torsion if its Witt equivalence class is a torsion element in the Witt group of quadratic forms over $K$.)

In this paper we  study these two invariants in the case of a Laurent series field $k(\!(t_1,\dotsc, t_n)\!)$
in $n\ge 2$ variables over a field $k$. (The $n=1$ case is classical.) Significant results in this direction already appeared in an influential paper of Choi, Dai, Lam and Reznick \cite{CDLR}. We exploit some newly developed methods to prove generalizations of some results in that paper.

\vskip2mm

The first main result is the following:

\begin{thm}
  \label{thm1p1hu}
Let $k$ be a field of characteristic different from $2$. Then
\[
p(k(\!(x,y)\!))=\sup\{p(\ell(x))\,|\,\;\ell/k \;\text{ a finite field extension }\}\,
\]and
\[
u(k(\!(x,y)\!))=2\sup\{u(\ell(x))\,|\,\;\ell/k \;\text{ a finite field extension }\}\,.
\]
\end{thm}
The proof will be completed in section$\;$\ref{secSeries2}. The starting point is a result in \cite{CDLR} which compares the sums of squares in $k(\!(x, y)\!)$ and those in $k(\!(y)\!)(x)$. As direct consequences, we get the inequalities
\[
p(k(\!(x,y)\!))\le p(k(\!(y)\!)(x))\quad \text{and }\quad u(k(\!(x,y)\!))\le u(k(\!(y)\!)(x))\,.
\]The equalities
\[
p(k(\!(y)\!)(x))=\sup\{p(\ell(x))\,|\,\;\ell/k \;\text{ a finite field extension }\}
\]and
\[
u(k(\!(y)\!)(x))=2\sup\{p(\ell(x))\,|\,\;\ell/k \;\text{ a finite field extension }\}\,
\]have been obtained recently by Becher, Grimm and Van Geel \cite{BGvG12}, using a local-global principle proved by Colliot-Th\'el\`ene, Parimala and Suresh \cite{CTPaSu} and some valuation-theoretic arguments.
These imply the  ``$\le$'' part of Theorem$\;$\ref{thm1p1hu}.
We will prove the inequalities in the other direction
by showing that each $\ell(x)$ is the residue field of a discrete valuation on $k(\!(x,y)\!)$ (see  Lemma$\;$\ref{lemma2p2hu}).

The statement on the Pythagoras number in Theorem$\;$\ref{thm1p1hu} generalizes the equality $p(\mathbb{R}(\!(x,y)\!))=2$ and the inequality $p(\mathbb{Q}(\!(x,y)\!))\le 8$, which were shown in \cite[$\S$5]{CDLR}. Our theorem implies that $p(\mathbb{Q}(\!(x,y)\!))=5$, since $p(\mathbb{Q}(x))=5$ and $p(\ell(x))\le 5$ for every finite extension $\ell$ of $\mathbb{Q}$ (cf. \cite[Chap.$\;$7, Thm.$\;$1.9]{Pfister95}).  This seems to give the first example of a formally real field $k$ for which $p(k(\!(x,y)\!))$ is not a power of 2.

The result on the $u$-invariant in Theorem$\;$\ref{thm1p1hu} can be viewed as a generalization of the equality $u(\mathbb{C}(\!(x,y)\!))=4$, first proved in \cite[Thm.$\;$5.16]{CDLR}. The $u$-invariant of two-variable Laurent series fields $k(\!(x,y)\!)$ and their finite extensions has been studied in a number of recent papers, e.g., \cite{CTOP}, \cite{HHK}, \cite{HHK11b}, \cite{Le10} and \cite{Hu11}. Most of these results deal with the case where $k$ is a nonreal field (i.e. a field in which $-1$ is a sum of squares), and our formula is new already in this case.

\vskip2mm

As another result on the Pythagoras number, we also give a generalization of
 the equality $p(\mathbb{R}(\!(x,y)\!))=2$ in a different direction. Namely, we show in Theorem$\;$\ref{thm5p1hu} that
 $p(L)\le 3$ for every finite extension $L$ of $\mathbb{R}(\!(x,y)\!)$. This is actually equivalent to the second equality in the following theorem.

\begin{thm}
  \label{thm1p2hu}
$p(\mathbb{R}(\!(x,y,z)\!))=p(\mathbb{R}(\!(x,y)\!)(z))=4$.
\end{thm}

Choi, Dai, Lam and Reznick  conjectured that $p(\mathbb{R}(\!(t_1,\dotsc, t_n)\!))\le 2^{n-1}$ for every $n\ge 3$ (cf. \cite[p.80, $\S$9, Problem$\;$6]{CDLR}). The best upper bound until now (even for $n=3$) is $2^n$, as shall be explained in section$\;$\ref{secSeries3}.

Theorem$\;$\ref{thm1p2hu}  is proved using a local-global principle   for isotropy of quadratic forms in 3 or 4 variables over finite extensions of $\mathbb{R}(\!(x,y)\!)$ (cf. \cite{Hu10}). We remark that  the statement of Theorem\;\ref{thm1p2hu} is still true when $\mathbb{R}$ is replaced by an iterated
Laurent series field $k=\mathbb{R}(\!(x_1)\!)\cdots (\!(x_m)\!)$ (Corollary$\;$\ref{coro5p3hu}).

 At the end of the paper, we propose two conjectures which predict that the formulas in Theorem$\;$\ref{thm1p1hu} have analogs for Laurent series fields in three or more variables.

\section{Lower bounds using discrete valuations}

Let $K$ be a field of characteristic $\neq 2$. A discrete valuation $v$ of  $K$ is called \emph{nondyadic}  if the residue field $\kappa(v)$ of $v$ has characteristic different from $2$.
It has been noticed by several authors that the invariants $p(K)$ and $u(K)$ can be bounded from below
in terms of those of the residue fields $\kappa(v)$, $v$ ranging over the nondyadic discrete valuations of $K$.

Unless otherwise stated,  we will follow standard notation for quadratic forms as used in \cite{Lam}.
As in \cite{BGvG12}, to avoid case distinction in some statements we set
\[
p'(K)=\begin{cases}
p(K)\;&\text{ if $K$ is (formally) real}\;\\
s(K)+1\;& \text{ if $K$ is nonreal}\,
\end{cases}
\]where $s(K)$ denotes the level of $K$ (see e.g. \cite[$\S$IX.2]{Lam}).

\begin{prop}[{See \cite[Propositions$\;$4.3 and 5.2]{BGvG12}, \cite[Prop.$\;$5]{Sch09}}]\label{prop2p1hu}
Let $v$ be a nondyadic discrete valuation of a field $K$.  Then
\[
p'(K)\ge p(K)\ge p'(\kappa(v))\quad\text{and }\quad u(K)\ge 2.u(\kappa(v))\,.
\]The equalities hold if $v$  is henselian $($meaning that the discrete valuation ring associated to $v$ is henselian$)$.
\end{prop}

This proposition generalizes the classical facts
\[
p'(K(\!(t)\!))=p(K(\!(t)\!))=p'(K)\,\quad\text{and }\quad u(K(\!(t)\!))=2u(K)
\](cf. \cite[Examples$\;$XI.5.9 (6) and Remarks$\;$XI.6.28 (2)]{Lam}).

Also, if $L/K$ is a finite separable field extension, then there is a discrete valuation $v$ on the rational function field $K(x)$ whose residue field $\kappa(v)$ is isomorphic to $L$. So, by Proposition$\;$\ref{prop2p1hu}, $p(K(x))\ge p'(L)$ and $u(K(x))\ge 2 u(L)$. The same is true when $L/K$ is a finite purely inseparable extension. For the Pythagoras number, this follows simply because $L$ is nonreal in this case and thus
\[
p'(L)=s(L)+1\le s(K)+1=s(K(x))+1=p(K(x))\,.
\]The argument for the $u$-invariant was given in \cite[Corollary$\;$5.4]{BGvG12}. Hence,  we have
\[
p(K(x))\ge \sup\{p'(L)\,|\, L/K \; \text{ a finite field extension }\}
\]and
\[
u(K(x))\ge 2\sup\{u(L)\,|\, L/K \; \text{ a finite field extension }\}\,.
\]These inequalities were shown in \cite{BGvG12} to be equalities when $K$ is a Laurent series field in one variable; this fact will be used in the next section.

We shall now consider the fraction field $K$ of a regular local ring with residue field $k$ and show that some special algebraic function fields over $k$ arise as residue fields of discrete valuations on $K$.

\begin{lemma}\label{lemma2p2hu}
  Let $A$ be a regular local ring of Krull dimension $n\ge 2$, and let $K$ and $k$ be respectively the fraction field and the residue field of $A$. Let $\ell/k$ be a finite field extension.

  Then the rational function field $\ell(x_1,\dotsc, x_{n-1})$ is the residue field of a discrete valuation on $K$.
\end{lemma}
\begin{proof}The proof makes use of the following fact: If $L/F$ is a finite simple extension of fields, then the affine line $\mathbb{A}^1_F$ has a closed point with residue field $L$. The same is true for any algebraic $F$-variety that contains $\mathbb{A}^1_F$ as a (locally closed) subvariety.

Now we use a geometric construction to derive a discrete valuation on $K$ with the given residue field. Let $X=\mathrm{Spec}(A)$ and let $X'\to X$ be the blowup of $X$ at its closed point. The exceptional divisor $E$ in $X'$ is isomorphic to $\mathbb{P}^{n-1}_k$ (cf. \cite[Thm.$\;$8.1.19]{Liu}). Let $k=\ell_0\subseteq \ell_1\subseteq \cdots\subseteq \ell_r=\ell$ be a chain of subfields of $\ell$ such that $\ell_{i+1}/\ell_i$ is a simple extension for each $i\in \{0,1,\dotsc, r-1\}$.

Since $n\ge 2$, there is a closed point  $Q\in E$ whose residue field $\kappa(Q)$ is isomorphic to $\ell_1$. In the blowup $X''$ of $X'$ at the point $Q$, we have an exceptional divisor $E'$ which is isomorphic to $\mathbb{P}^{n-1}_{\ell_1}$. We can choose a closed point on $E'$ whose residue field is $\ell_2$ and blow up $X''$ at that point. Then we get a regular scheme $X^{(3)}$ which is birational to $X$ and which contains a divisor  isomorphic to $\mathbb{P}^{n-1}_{\ell_2}$. Repeating this procedure sufficiently many times produces a regular scheme  birational to $X$ which contains $\mathbb{P}^{n-1}_{\ell}$ as a divisor. The generic point of this divisor defines a discrete valuation of $K$ whose residue field is $\ell(x_1,\dotsc, x_{n-1})$.
\end{proof}
A variant of the above lemma has been noticed independently in \cite[Lemma$\;$4.1]{Grimm}, where the blowup construction is expressed in a purely algebraic form. In geometric terms, the proof there is based on the observation that for any simple extension $L/k$, the  blowup $X'$ of $\mathrm{Spec}(A)$ (at its closed point) has a point with residue field $L$.
Our proof of Lemma$\;$\ref{lemma2p2hu} uses this fact when $L/k$ is an algebraic simple extension.

By carrying out blowups over polynomials rings over $A$, it is possible to get similar results for discrete valuations on rational function fields over $K$.

Applying Proposition$\;$\ref{prop2p1hu} and Lemma$\;$\ref{lemma2p2hu} to the power series ring $k[\![t_1,\dotsc, t_n]\!]$ and its fraction field $k(\!(t_1,\dotsc, t_n)\!)$, we obtain the following corollary.

\begin{coro}\label{coro2p3hu}
For any field $k$ of characteristic different from $2$ and any $n\ge 2$, one has
\[
p(k(\!(t_1,\dotsc, t_n)\!))\ge \sup\{p(\ell(x_1,\dotsc, x_{n-1}))\,|\, \ell/k \; \text{ a finite field extension}\}\,
\]and
\[
u(k(\!(t_1,\dotsc, t_n)\!))\ge 2\sup\{u(\ell(x_1,\dotsc, x_{n-1}))\,|\,\ell/k \; \text{ a finite field extension}\}\,.
\]
\end{coro}
In \cite[Theorem$\;$3.3]{Grimm} the same lower bound for the Pythagoras number was shown for algebraic function fields of transcendence degree $n$ over $k$ in place of Laurent series fields in $n$ variables.

\section{Laurent series in two variables}\label{secSeries2}

The goal of this section is to prove Theorem$\;$\ref{thm1p1hu} and to give some applications.
The following analogous result will be used in the proof of our theorem.

\begin{thm}\label{thm3p1hu} Let $k$ be a field of characteristic different from $2$.

$(\oi)$ $p(k(x))\le p(k(\!(t)\!)(x))$ and these two Pythagoras numbers are bounded by the same $2$-powers.

$(\ii)$ $p(k(\!(t)\!)(x))=\sup\set{p(\ell(x))\,|\, \ell/k \text{ a finite field extension}}$.

$(\iii)$ $u(k(\!(t)\!)(x))=2\sup\set{u(\ell(x))\,|\, \ell/k \text{ a finite field extension}}$.
\end{thm}
\begin{proof}
(i) The first assertion is contained in \cite[Prop.$\;$5.17]{Sch01}. The proof for the second assertion already appeared in \cite[Thm.$\;$5.18]{CDLR}. (See also \cite[Thm.$\;$4.14]{BGvG12}.)

(ii)  \cite[Coro.$\;$6.9]{BGvG12}.

(iii) \cite[Thm.$\;$6.6]{BGvG12}.
\end{proof}

It is conjectured in \cite[Conjecture$\;$4.16]{BGvG12} that the inequality in assertion (i) of Theorem$\;$\ref{thm3p1hu} is actually an equality or, equivalently, that $p(\ell(x))\le p(k(x))$ for all finite extensions $\ell/k$.

Given a field $k$,  we will frequently write $R_n$ for the ring of formal power series $k[\![t_1,\dotsc, t_n]\!]$ in the variables $t_1,\dotsc, t_n$ for each $n\ge 1$. $F_n$ will be a shorthand for $k(\!(t_1,\dotsc, t_n)\!)$, the corresponding field of Laurent series. By convention $R_0=F_0=k$.

We start the proof of Theorem$\;$\ref{thm1p1hu} with the assertion on the Pythagoras number in the nonreal case. We can prove the following more general fact.

\begin{prop}\label{prop3p2hu}
If $k$ is a nonreal field of characteristic different from $2$, then for every $n\ge 1$ one has
\[
s(F_n)=s(k)\quad\text{ and } \quad p(F_n)=s(k)+1=p(F_n(t))\,.
\]
\end{prop}
\begin{proof}
Since $p(K(t))=s(K)+1$ for any nonreal field $K$,  we need only prove the  equalities $s(F_n)=s(k)$ and $p(F_n)=s(k)+1$. We use induction on $n$. The case $n=1$, as discussed previously, is a special case of Proposition$\;$\ref{prop2p1hu}.

Assume $n\ge 2$. The inclusions $k\subseteq F_n\subseteq F_{n-1}(\!(t_n)\!)$ yield
\[
s(F_{n-1}(\!(t_n)\!))\le s(F_n)\le s(k)\,.
\]But $s(F_{n-1}(\!(t_n)\!))=s(F_{n-1})=s(k)$ by the $n=1$ case and the induction hypothesis. This proves $s(F_n)=s(k)$.

For the Pythagoras number, we have  $p(F_n)\le s(F_n)+1=s(k)+1$. On the other hand, Lemma$\;$\ref{lemma2p2hu} together with Proposition$\;$\ref{prop2p1hu} shows that
\[
p(F_n)\ge p'(k(x_1,\dotsc, x_{n-1}))=s(k(x_1,\dotsc, x_{n-1}))+1=s(k)+1\,.
\]Alternatively, one can prove the inequality $p(F_n)\ge s(k)+1$ by showing that
$t_n$ cannot be expressed as a sum of $s(k)$ squares in the field $F_{n-1}(\!(t_n)\!)$.
\end{proof}

For a nonreal field $k$, one has
\[
p(k(x))=s(k)+1\ge s(\ell)+1=p(\ell(x))
\]for any finite extension $\ell/k$. So the result on the Pythagoras number in Theorem$\;$\ref{thm1p1hu} in the nonreal case is covered by Proposition$\;$\ref{prop3p2hu}.
The real case will be treated in Theorem$\;$\ref{thm3p4hu}.

\begin{lemma}\label{lemma3p3hu}
Let $k$ be a field of characteristic different from $2$. Consider the rings $k[\![t]\!][x]\subseteq k[x][\![t]\!]\subseteq k[\![x,t]\!]$.

$(\mathrm{i})$ Every  $f\in k[x][\![t]\!]\,($resp. $f\in k[\![x,t]\!])$ admits a factorization $f=u^2g$, where $u$ is a unit in $k[x][\![t]\!]\,($resp. in $k[\![x,t]\!])$ and $g\in k[\![t]\!][x]$.

$(\mathrm{ii})$ Suppose $k$ is real. Then for every $m\ge 1$, every sum of $m$ squares in $k[x][\![t]\!]\,($resp. in $k[\![x,t]\!])$ is of the form $a^2b$, where $a$ lies in $k[x][\![t]\!]($resp. in $k[\![x,t]\!])$ and $b$ is a sum of $m$ squares in $k[\![t]\!][x]$.
\end{lemma}
\begin{proof}
(i) In the ring $k[x][\![t]\!]$ (resp. $k[\![x,t]\!]$) every unit is the product of a square and an element in $k^*\subseteq k[\![t]\!][x]$. So it is enough to factorize $f$ as $f=ug$ with $u$ a unit in $k[x][\![t]\!]$ (resp. $k[\![x,t]\!]$)  and $g\in k[\![t]\!][x]$.
We may assume $t\nmid f$.

Then the statement for the ring $k[\![x,t]\!]$ is classical (cf. \cite[p.145, Coro.$\;$1]{ZS2}).

  Let us consider the statement for the ring $k[x][\![t]\!]$. Write $f=\sum_{i\ge 0}f_i(x)t^i$ with $f_i\in k[x]$.  Applying a Weierstrass lemma as stated in \cite[Lemma$\;$5.3]{CDLR} to the rings
  $A=k[x][\![t]\!]$, $B=k[x]$ and the subspace $C\subseteq B$ of polynomials of degree $<d:=\deg(f_0)$ together with $p=f$, we get  an expression
$x^d=qf+r$, where $q\in k[x][\![t]\!]$ and $r\in k[\![t]\!][x]$ has degree $<d$ in $x$. Considering this equation modulo $t$, we see that the constant term $q_0\in k[x]$ of $q\in k[x][\![t]\!]$ is of degree 0 in $x$. Hence, $q$ is a unit in $k[x][\![t]\!]$ and we can take $u=1/q$ and $g=x^d-r$.

(ii) See \cite[Thm.$\;$5.20]{CDLR}.
\end{proof}

The Pythagoras number of a  ring can be defined as in the case of a field (cf. \cite[p.45]{CDLR}).
If $A$ is an integral domain, we denote by $\mathrm{Frac}(A)$ its fraction field.

The following result strengthens \cite[Coro.$\;$5.21]{CDLR}.

\begin{thm}\label{thm3p4hu}
Let $k$ be a real field. Then the rings
\[
k[\![x,\,t]\!]\,,\;  k[x][\![t]\!]\,,\;k[\![t]\!][x]\,,\; k(\!(x,\,t)\!)\,,\; \mathrm{Frac}(k[x][\![t]\!])\;\text{ and }\;\; k(\!(t)\!)(x)
\]have the same Pythagoras number, which is equal to
\[
 \sup\{p(\ell(x))\,|\,\ell/k \;\text{ a finite field extension }\}\,.
 \]
\end{thm}
\begin{proof}
From Lemma$\;$\ref{lemma3p3hu} (ii) it follows that
\[
p(k(\!(x,\,t)\!))\le p(k[\![x,\,t]\!])\le p(k[x][\![t]\!]) \le p(k[\![t]\!][x])\]
 and
\[
p(k(\!(x,\,t)\!))\le p(\mathrm{Frac}(k[x][\![t]\!]))
\le p(k(\!(t)\!)(x))\,.
\]On the other hand, the proof of \cite[Thm.$\;$5.18]{CDLR} has actually shown the equality $p(k[\![t]\!][x])=p(k(\!(t)\!)(x))$.
So it suffices to show
\[
p(k(\!(x,\,t)\!))\ge \sup\{p(\ell(x))\,|\,\ell/k \;\text{ a finite field extension }\}=p(k(\!(t)\!)(x))\,.
\]But this follows by combining Theorem$\;$\ref{thm3p1hu} (ii) and Corollary$\;$\ref{coro2p3hu}.
\end{proof}

\begin{remark}
Recall that the earlier mentioned \cite[Conjecture$\;$4.16]{BGvG12} would imply (equivalently) that
one can replace $\sup\{p(\ell(x))\,|\,\ell/k\, \text{ finite extension}\}$ by $p(k(x))$ in the statement of Theorem$\;$\ref{thm3p4hu}, i.e. that one would obtain the equality $p(k(\!(x,y)\!))=p(k(x))$. We observe that this is indeed the case when $k$ is either real closed (e.g. $k=\mathbb{R}$) or a number field (i.e. finite extension of $\mathbb{Q}$).

The equality $\sup\{p(\ell(x))\,|\,\ell/k\, \text{ finite extension}\}=p(k(x))$ is trivially true when $k$ is real closed. For $k$ a nonreal number field it follows from Proposition$\;$\ref{prop3p2hu}, and the case where $k$ is a real number field follows either from Theorem$\;$\ref{thm3p1hu} (i) (when $p(k(x))=4$), or directly from the general result  $4\le p(\ell(x))\le 5$ for any real number field $\ell$, due to Pourchet and Hsia--Johnson (cf. \cite[Chap.$\;$7, Thm.$\;$1.9]{Pfister95}). In particular, we have
\[
p(\mathbb{Q}(\!(x,y)\!))=5\,.
\]
\end{remark}

For a field $K$ of characteristic different from 2, we denote by $(\sum K^2)^*$ the multiplicative group of nonzero sums of squares in $K$. This is the same as the group of totally positive elements in $K^*$ by \cite[Thm.$\;$VIII.1.12]{Lam}, and coincides with $K^*$ if $K$ is nonreal.

The author thanks K.\! Becher for helpful discussions on the following lemma, which should be well known to experts.

\begin{lemma}\label{lemma3p6hu}\footnote{The present form of this lemma is kindly suggested by an anonymous referee.}
  Let $L/K$ be an extension of fields of characteristic different from $2$ such that the  natural homomorphism
  $\varphi: K^*/K^{*2} \to L^*/L^{*2}$ is surjective.

  Then the natural homomorphism of Witt rings $\psi: W(K)\to W(L)$ is surjective.

  If furthermore the restriction of $\varphi$ to totally positive squares classes $(\sum K^2)^*/K^{*2}\to (\sum L^2)^*/L^{*2}$ is surjective, then so is the restriction of $\psi$ to the torsion parts of the fundamental ideals of Witt classes of even dimensional forms $I(K)_{tors}\to I(L)_{tors}$.
\end{lemma}
\begin{proof}
The Witt group of a field is generated by one-dimensional forms, so the first assertion follows immediately. For the second assertion, we use the fact that the torsion part of the fundamental ideal is generated as a group by two-dimensional forms of the shape $\alpha.\dgf{1,\,-\beta}=\dgf{\alpha,\,-\alpha\beta}$, where $\alpha$ is an arbitrary nonzero field element and $\beta$ is a totally positive one. The latter fact is trivially true for nonreal fields as every two-dimensional form is of that shape, and it is shown in
\cite[Satz.$\;$22]{Pfister66} in the case of real fields.
\end{proof}

We can now prove the following theorem, which covers the result on the $u$-invariant in Theorem$\;$\ref{thm1p1hu}.

\begin{thm}\label{thm3p7hu}
For any  field $k\,$of characteristic different from $2$, one has
\[
u(k(\!(x,\,t)\!))=u(\mathrm{Frac}(k[x][\![t]\!]))=u(k(\!(t)\!)(x))=2.\sup\{u(\ell(x))\,|\,\ell/k \text{ a finite extension}\}\,.
\]
\end{thm}
\begin{proof}
The last equality was proved in \cite[Thm.$\;$6.6]{BGvG12}.

We first show the inequalities $u(k(\!(x,\,t)\!))\le u(\mathrm{Frac}(k[x][\![t]\!]))\le u(k(\!(t)\!)(x))$. It is sufficient to prove that the corresponding homomorphisms given by scalar extension
\[
W(k(\!(t)\!)(x))_{tors}\longrightarrow W(\mathrm{Frac}(k[x][\![t]\!]))_{tors}\longrightarrow W(k(\!(x,\,t)\!))_{tors}\,.
\]are surjective. Note that all these fields are real if and only if $k$ is real. In the nonreal case, the Witt groups are torsion groups (cf. \cite[Corollary$\;$XI.2.3]{Lam}), and in the real case, the torsion parts of the Witt groups coincide with the torsion parts of the fundamental ideals as can be easily seen from \cite[Thm.$\;$VIII.3.2]{Lam}. That  the first homomorphism and the composition (and hence also the second homomorphism) are surjective follows from Lemma$\;$\ref{lemma3p6hu}, as soon as we verified that the hypotheses on the natural maps of square classes are satisfied for the respective fields. The latter is done in Lemma$\;$\ref{lemma3p3hu}.

Finally, we have
\[
u(k(\!(t)\!)(x))=2.\sup\{u(\ell(x))\,|\,\ell/k \text{ a finite extension}\}\le u(k(\!(x,\,t)\!))
\]by Corollary$\;$\ref{coro2p3hu}. This completes the proof.
\end{proof}

\begin{remark}\label{remark3p8hu}
(1) If $k$ is a nonreal field, the inequality $u(k(\!(x,\,t)\!))\le u(k(\!(t)\!)(x))$ is implicitly contained in \cite{CDLR}.
However, even in the nonreal case the first two equalities in Theorem$\;$\ref{thm3p7hu} seem to have escaped earlier notice. Moreover,
the relation $u(k(\!(x,\,t)\!))=u(k(\!(t)\!)(x))$ has a mixed characteristic version: For a complete discrete valuation ring $A$, one has
\[
u(\mathrm{Frac}(A[\![x]\!]))=u(\mathrm{Frac}(A[x]))\,,
\]as can be shown analogously (using \cite[Thm.$\;$6.6]{BGvG12} and an analog of Corollary$\;$\ref{coro2p3hu}).

(2) If $u_s(k)$ denotes the \emph{strong $u$-invariant} of $k$ as defined in  \cite[Definition$\;$2]{Sch09} or \cite[$\S$5]{BGvG12}, i.e.,
\[
u_s(k)=\frac{1}{2}\sup\{u(L)\,|\,L/k \text{ a finitely generated extension of transcendence degree }1 \}\,,
\]
then Theorem$\;$\ref{thm3p7hu} implies that $u(k(\!(x,\,t)\!))\le 4u_s(k)$. In the nonreal case, \cite[Coro.$\;$4.2]{HHK11b} gives a generalization of this inequality.
\end{remark}

Let us close this section with some examples where our results, combined with some recent work of others, can give refinements of earlier results (especially those obtained in \cite[$\S$5]{CDLR}).

\begin{example}\label{exam3p9hu}
   Let $k_0$ be a real closed field or a number field, and let $k$ be a finitely generated extension of transcendence degree $d\ge 1$ over $k_0$. Then
   \[
\sup\{p(\ell(x))\,|\,\ell/k \text{ a finite extension}\}\le \begin{cases}
  2^{d+1}\; & \text{ if $k_0$ is real closed},\\
  2^{d+2}\; & \text{ if $k_0$ is a number field}.
\end{cases}
\]
   When $k_0$ is real closed, the upper bound is due to Pfister (cf. \cite[Chap.$\;$7, Examples$\;$1.4 (4)]{Pfister95}).   When $k$ is a number field, it was shown in \cite[Thm.$\;$4.1]{CTJannsen91} assuming Milnor's conjecture and another conjecture by Kato that have both been solved in \cite{OVV07} and \cite{Jannsen09} respectively.

On the other hand, if $k$  is real, one has (see e.g. \cite{Grimm})
\[
p(k(x))\ge \begin{cases}
  d+2\; & \text{ if $k_0$ is real closed},\\
  d+3\; & \text{ if $k_0$ is a number field}.
\end{cases}
\](By \cite[Prop.$\;$2.2]{Grimm}, there is a finite real extension $k'_0$ of $k_0$ such that
$k_0'(x_1,\dotsc, x_d)$ arises as the residue field of a discrete valuation on $k(x)$. In the number field case, this implies that $p(k(x))\ge p(k'_0(x))+d-1\ge d+3$.)

Therefore, for $k$ real, one has
\[
\begin{cases}
  d+2\le p(k(\!(x,y)\!))\le 2^{d+1}\; & \text{ if $k_0$ is real closed},\\
  d+3\le p(k(\!(x,y)\!))\le 2^{d+2}\; & \text{ if $k_0$ is a number field}.
\end{cases}
\]
In some special cases, one  can get a better lower bound. For example, if $k=\mathbb{Q}(x_1,\dotsc, x_d)$, then
$p(k(\!(x,y)\!))\ge p(k(x))\ge d+p(\mathbb{Q}(x))=d+5$.

\vskip2mm

In the case where $k_0$ is real closed, one has
\[
2^{d+1}\le u(L)\le 2^{d+3}-2d-8\,
\] for every finitely generated extension $L/k_0$ of transcendence degree $d+1$. The upper bound was shown in \cite[Thm.$\;$4.11]{ElmanLam73MZ} if $d=1$ and follows from  \cite[Thm.$\;$3]{Becher10} if $d\ge 2$. Thus,  Theorem$\;$\ref{thm3p7hu} (or Theorem$\;$\ref{thm1p1hu}) yields
\[
2^{d+2}\le u(k(\!(x,y)\!))\le 2^{d+4}-4d-16\,.
\]In particular,
\[
8\le u(\mathbb{R}(t)(\!(x,\,y)\!))\le 12\,.
\]The upper bound here is sharper than the upper bound one can get directly from \cite[Thm.$\;$3]{Becher10} (without using Theorem$\;$\ref{thm3p7hu}).
\end{example}

\begin{example}
This example gives a new class of fields $k$ for which $p(k(\!(x,y)\!))=4$.

Let $k_0$ be a real closed field. For each $m\ge 0$, let $k_m=k_0(\!(x_1)\!)\cdots (\!(x_m)\!)$ be an iterated Laurent series field over $k_0$.
We will show in Corollary$\;$\ref{coro5p3hu}
that $p(k_m(\!(t_1,\,t_2)\!)(x))=4$ for every $m$. Combining Theorem$\;$\ref{thm3p1hu} (i) and Theorem$\;$\ref{thm3p4hu}, we get $p(k_m(\!(t_1,\,t_2)\!)(\!(x,y)\!))=4$.

In particular,
\[
p(\mathbb{R}(\!(x_1)\!)\cdots (\!(x_m)\!)(\!(t_1,\, t_2)\!)(\!(x,\,y)\!))=4\,,\quad \forall\; m\ge 0\,.
\]
\end{example}

\section{Comparison results using Weierstrass preparation}\label{Weierstrass}

In this section,  we  use Weierstrass-type theorems for power series in several variables to generalize some comparison results in the previous section.

Let $k$ be a field  and $n\ge 1$ an integer. A power series $f\in R_n=k[\![t_1,\dotsc, t_n]\!]$ is said to be \emph{regular in $t_n$} if $f(0,\dotsc, 0,\,t_n)\neq 0$ in $k[\![t_n]\!]$.

\begin{lemma}\label{lemma4p1hu}Assume the field $k$ has characteristic different from $2$.
Let $f\in R_n$ be a power series regular in $t_n$ and let $m$ be an integer such that $1\le m\le s(k)$.

If $f$ is a sum of $m$ squares in $R_n$, then there is a factorization $f=a^2b$, where
 $a \in R_n$ and  $b$ is a sum of $m$ squares in $R_{n-1}[t_n]$.
\end{lemma}
\begin{proof}The proof is similar to that of \cite[Thm.$\;$5.20]{CDLR}.

Write $f=\sum_{i=1}^mg_i^2$ with $g_i\in R_n$. Let $d$ be the $t_n$-adic valuation of $f(0,\dotsc, 0,\,t_n)\in k[\![t_n]\!]$. If $d=0$, then $f$ is a unit in $R_n$ and
$\alpha:=f(0,\dotsc, 0)$ is a nonzero element in $k$. In this case  $\alpha^{-1}f$ is a square in $R_n$, and clearly $\alpha=\sum^m_{i=1}g_i(0,\dotsc, 0)^2$ is a sum of $m$ squares in $R_{n-1}[t_n]$. So we may take $b=\alpha$ and $a\in R_n$ a square root of $\alpha^{-1}f$.

Now assume $d>0$, so that $0=\sum^m_{i=1}g_i(0,\dotsc, 0)^2$ in $k$. Since $m\le s(k)$ by assumption, one has $g_i(0,\dotsc, 0)=0$ for every $i$. By \cite[p.139, Thm.$\;$5]{ZS2},  each $g_i$ can be  written as $g_i=h_if+r_i$, where $h_i\in R_n$ and $r_i\in R_{n-1}[t_n]$ with $\deg_{t_n}(r_i)<d$. Then the power series
\[
u:=1-2\sum_{i=1}^mg_ih_i+f\sum_{i=1}^mh_i^2\,,
\] is a unit and is a square in $R_n$ since $u(0,\dotsc, 0)=1$, and one has
\[
\sum r_i^2=\sum(g_i-h_if)^2=\sum g_i^2-2f\sum g_ih_i+f^2\sum h_i^2=fu\,.
\]Now taking $a\in R_n$ such that $a^2=u^{-1}$ and $b=\sum r_i^2$ finishes the proof.
\end{proof}

For a commutative ring $A$ and an integer $m\ge 1$, we denote by $D_A(m)$ the set of nonzero sums of $m$ squares in $A$.
\begin{prop}\label{prop4p2hu}
Let $k$ be a real field and $m\ge 1$ an integer.

Then given finitely many elements $f_1,\dotsc, f_r\in D_{R_n}(m)$, there is an automorphism $\sigma$ of the ring $R_n$ such that $\sigma(f_i)\in R^2_n\cdot D_{R_{n-1}[t_n]}(m)$ for every $1\le i\le r$. If $n\le 2$, we may take $\sigma$ to be the identity.
\end{prop}
\begin{proof}
The case with $n=1$ is left to the reader. If $n=2$, this is part of  Lemma$\;$\ref{lemma3p3hu}. For general $n$, it suffices to apply \cite[p.147, Corollary]{ZS2} to get an automorphism $\sigma$ such that all the $\sigma(f_i)$ are regular in $t_n$. Then the result follows from Lemma$\;$\ref{lemma4p1hu}.
\end{proof}

In general, the automorphism $\sigma$ in the above proposition may not preserve the subring $R_{n-1}[t_n]$.

\begin{prop}\label{prop4p3hu}Assume $k$ is of characteristic different from $2$.
For every torsion form  $\phi$ over $F_n=k(\!(t_1,\dotsc, t_n)\!)$, there is an automorphism $\sigma$ of $F_n$ and a torsion form $\psi$ over $F_{n-1}(t_n)$ such that
$\sigma_*[\phi]=[\psi]$ in $W(F_n)$, where
\[
\sigma_*\,:\;W(F_n)\longrightarrow W(F_n)\,;\quad\dgf{f_1\,,\dotsc, f_r}\longmapsto \dgf{\sigma(f_1)\,,\dotsc, \sigma(f_r)}
\]denotes the automorphism of the Witt group $W(F_n)$ induced by $\sigma$. If $n\le 2$, one can take $\sigma$ to be the identity.
\end{prop}
\begin{proof}
We may assume $\phi=\dgf{f_1,\dotsc, f_r}$, where all the coefficients $f_i$ lie in $R_n$.

First assume $k$ is a nonreal field. By \cite[p.145, Coro.$\;$1 and p.147, Corollary]{ZS2}, there is an automorphism $\sigma$ of $R_n$ such that for every $i$ the power series $\sigma(f_i)$ admits a factorization $\sigma(f_i)=u_ig_i$, where $u_i$ is a unit in $R_n$ and $g_i\in R_{n-1}[t_n]$. When $n\le 2$, one can take $\sigma$ to be the identity. (If $n=2$, one has $f_i=t_1^{r_i}f'_i$ for some $r_i\ge 0$ and $f'_i$ regular in $t_2$.)
 Putting
\[
\lambda_i:=u_i(0,\dotsc, 0)\in k^*\quad\text{ and }\quad h_i:=\lambda_ig_i\,,
\]we get a form $\psi:=\dgf{h_1,\dotsc, h_r}$ which is defined over $F_{n-1}(t_n)$ and isomorphic to
\[
\sigma(\phi):=\dgf{\sigma(f_1)\,\dotsc, \sigma(f_r)}\,.
\] over $F_n$.

Now consider the case with $k$ real. Then the form $\phi$ is Witt equivalent to  $\psi_1\bot\cdots\bot\psi_m$ for some binary torsion forms $\psi_j=c_j.\dgf{1,\,-d_j}$ with $c_j,\,d_j\in R_n$. Each $d_j$ is a sum of squares in $F_n$, and we may assume it is already a sum of squares in $R_n$.

By Proposition$\;$\ref{prop4p2hu}, there is an automorphism $\sigma$ of $R_n$, which we may take to be the identity if $n\le 2$, such that each of the $\sigma(d_j)$ is a sum of squares in $R_{n-1}[t_n]$ up to a square in $R_n$. We may and we will further assume that $\sigma$ is chosen such that each
$\sigma(c_j)$ admits a factorization $\sigma(c_j)=u_je_j$, where $u_j$ is a unit in $R_n$ and $e_j\in R_{n-1}[t_n]$. (For $n\ge 3$, it suffices to choose $\sigma$ such that all the $\sigma(d_j)$ and $\sigma(c_j)$ are regular in $t_n$.) Thus, there are elements $d_j',\,c_j'\in R_{n-1}[t_n]$ such that
$d_j'$ is a sum of squares in $R_{n-1}[t_n]$ and
\[
d_j'\sigma(d_j)^{-1}\in R_n^2\,,\quad c_j'\sigma(c_j)^{-1}\in R_n^2
\]
for every $j$. Now the form
\[
\psi:=c'_1.\dgf{1,\,-d'_1}\,\bot\cdots\bot\, c'_m.\dgf{1,\,-d'_m}
\]is a torsion form over $F_{n-1}(t_n)$ with the desired property.
\end{proof}

\begin{coro}\label{coro4p4hu}Assume $k$ is of characteristic different from $2$.
For every $n\ge 1$, one has  $p(F_n)\le p(F_{n-1}(x))$ and $u(F_n)\le u(F_{n-1}(x))$.
\end{coro}
\begin{proof}
First consider the assertion about the Pythagoras number. If $k$ is nonreal, then $p(F_n)=p(F_{n-1}(x))$ by Proposition$\;$\ref{prop3p2hu}. So we may assume $k$ is real and $m:=p(F_{n-1}(x))<\infty$. We want to show that every sum of squares $f$ in $F_n$ is a sum of $m$ squares. By Proposition$\;$\ref{prop4p2hu}, there is an automorphism $\sigma$ of $F_n$ such that $\sigma(f)\in D_{F_n}(m)$. This implies $f\in D_{F_n}(m)$.

For the assertion about the $u$-invariant, let $u=u(F_{n-1}(x))\le \infty$ and let $\phi$ be a torsion form over $F_n$. We need to show that the anisotropic part
$\phi_{an}$ of $\phi$ has dimension at most $u$. Indeed, Proposition$\;$\ref{prop4p3hu} provides an automorphism $\sigma$ of the field $F_n$ such that
$[\sigma(\phi)]=\sigma_*[\phi]=[\psi]\in W(F_n)$ for some torsion form $\psi$ over $F_{n-1}(t_n)$. Hence,
\[
\dim\phi_{an}=\dim\sigma(\phi)_{an}\le \dim\psi_{an}\le u=u(F_{n-1}(t_n))\,,
\] completing the proof.
\end{proof}

\section{Laurent series fields in three variables}\label{secSeries3}

  Pfister  proved that for any integer $r\ge 1$, if $F$ is a field of characteristic different from 2 such that every $(r+1)$-fold Pfister form over $F(\sqrt{-1})$ is hyperbolic, then $p(F)\le 2^r$ (cf. \cite[Chap.$\;$6, Thm.$\;$3.3]{Pfister95}). Combined with Milnor's conjecture (proved in \cite{OVV07}), this shows that if the  cohomological $2$-dimension $\mathrm{cd}_2(F(\sqrt{-1}))$ of the field $F(\sqrt{-1})$ is at most $r$, then $p(F)\le 2^r$. For example, for $F=\mathbb{R}(\!(t_1,\dotsc, t_n)\!)$, one has $\mathrm{cd}_2(F(\sqrt{-1}))\le n$ by \cite[XIX, Coro.$\;$6.3]{SGA4}. Hence
$p(\mathbb{R}(\!(t_1,\dotsc, t_n)\!))\le 2^n$ for all $n\ge 1$. In particular, this argument yields $p(\mathbb{R}(\!(t_1,\,t_2,\, t_3)\!))\le 8$.

Our goal here is to determine the precise value of $p(\mathbb{R}(\!(t_1,\,t_2,\, t_3)\!))$.

\begin{thm}\label{thm5p1hu}
Let $k$ be a real field such that $p(k(x,\,y))\le 4\,($e.g., $k$ real closed$)$. Then
\[
p(k(\!(t_1,\,t_2,\,t_3)\!))=p(k(\!(t_1,\,t_2)\!)(x))=4\,
\]and  $2\le p(L)\le 3$ for every finite extension $L$ of $k(\!(t_1,\,t_2)\!)$.
\end{thm}
\begin{proof}
In fact,  one has  $p(K(x,\,y))\ge 4$ for any real (but not necessarily real closed) field $K$. (It is classical that the Motzkin polynomial
(cf. \cite[Chap.$\;$1, Example$\;$2.4]{Pfister95}) is not a sum of fewer than 4 squares in $\mathbb{R}(x,\,y)$. It is not difficult to see that the same is true if $\mathbb{R}$ is replaced with an arbitrary real closed field. Since any real field $K$ is contained in a real closed field, the same result holds over $K(x,\,y)$.)  So the hypothesis $p(k(x,\,y))\le 4$ is equivalent to $p(k(x,\,y))=4$.

We have thus $p(k(\!(t_1,t_2,t_3)\!))\ge p(k(x,\,y))=4$ by Corollary$\;$\ref{coro2p3hu}.  (As pointed out by Becher, one can also prove $p(k(\!(t_1,t_2,t_3)\!))\ge 4$ by using the Motzkin polynomial.) In view of Corollary$\;$\ref{coro4p4hu}, the two equalities asserted in the theorem will follow from the inequality $p(k(\!(t_1,t_2)\!)(x))\le 4$. By \cite[Thm.$\;$3.5]{BvG09}, which is an extended version of a theorem of Pfister (cf. \cite[Examples$\;$XI.5.9 (3)]{Lam}),
this last condition is equivalent to each of the following two conditions:

(1) $p(L)\le 3$ for every finite extension $L$ of $k(\!(t_1, t_2)\!)$.

(2) $s(L)\le 2$ for every finite \emph{nonreal} extension $L$ of $k(\!(t_1, t_2)\!)$.

It is easy to see that finite extensions of $k(\!(t_1, t_2)\!)$ have Pythagoras number at least 2. So it remains to prove that condition (2) is satisfied.

Fix a finite nonreal extension $L/k(\!(t_1, t_2)\!)$ and let $R$ be the integral closure of $k[\![t_1,\,t_2]\!]$ in $L$. Call a discrete valuation $w$ of $L$ \emph{divisorial} if there is a regular integral scheme $X$ equipped with a proper birational morphism $X\to \mathrm{Spec}(R)$ such that $w$ is defined by a codimension 1 point of $X$. For such a discrete valuation $w$, the completion $L_w$ is nonreal as $L$ is, and
 the residue field $\kappa(w)$ is (isomorphic to) either a finite extension of $k(x)$ or the fraction field of a complete discrete valuation ring whose residue field $\ell$  is a finite extension of $k$. (This fact is probably evident to algebraic geometers. A detailed explanation can be found in \cite[Coro.$\;$2.3.26]{HuThese}.)
In the former case, the hypothesis on $k(x,\,y)=k(x)(y)$ implies that $s(\kappa(w))\le 2$ (by Pfister's theorem or \cite[Thm.$\;$3.5]{BvG09}). In the latter case, we have $s(\kappa(w))=s(\ell)\le 2$ since $p(k(x))\le p(k(x,\,y))\le 4$. Hence, we have $s(L_w)=s(\kappa(w))\le 2$  in any case.

It is proved in \cite[Thm.$\;$1.1]{Hu10} that the isotropy of quadratic forms of rank 3 or 4 over $L$ satisfies the local-global principle with respect to the divisorial valuations of $L$. This implies immediately $s(L)\le 2$ as desired.
\end{proof}

\begin{remark}\label{remark5p2hu}
(1) The hypothesis $p(k(x,\,y))\le 4$ in Theorem$\;$\ref{thm5p1hu} is satisfied if $k$ is a hereditarily euclidean field (cf. \cite[Coro.$\;$4.6]{BvG09}).

(2) If a field $k$ satisfies $p(k(x,\,y))\le 4$, then so does any iterated Laurent series field $k'=k(\!(x_1)\!)\cdots(\!(x_m)\!)$. Indeed, we may assume $k'=k(\!(x_1)\!)$ by induction and we need only to show $s(L)\le 2$ for every finite nonreal extension $L$ of $k'(x)$. The argument for a similar statement given in our proof of Theorem$\;$\ref{thm5p1hu} works verbatim, since we have a local-global principle for 3-dimensional quadratic forms over $L$ (see e.g. \cite[Thm.$\;$3.1]{CTPaSu}).
\end{remark}

\begin{coro}\label{coro5p3hu}
Let $k_0$ be a real closed field and $k=k_0(\!(x_1)\!)\cdots(\!(x_m)\!)$ an iterated Laurent series field over $k_0$.

Then we have
\[
p(k(\!(t_1,\,t_2,\,t_3)\!))=p(k(\!(t_1,\,t_2)\!)(x))=4\,
\]and  $2\le p(L)\le 3$ for every finite extension $L$ of $k(\!(t_1,\,t_2)\!)$.
\end{coro}
\begin{proof}
Immediate from Theorem$\;$\ref{thm5p1hu} and Remark$\;$\ref{remark5p2hu} (2).
\end{proof}

Note that for the field $k$ in the above corollary, Pfister's method (mentioned at the beginning of this section) can only give the upper bound
$2^{m+3}$.

\

Now consider the field $F_n=k(\!(t_1,\dotsc, t_n)\!)$ for general $n$.
For any $n\ge 2$, we have shown (in Corollaries$\;$\ref{coro2p3hu} and \ref{coro4p4hu})
\[
p(F_{n-1}(x))\ge p(F_n)\ge \sup\{p(\ell(x_1,\dotsc, x_{n-1}))\,|\,\ell/k \text{ a finite field extension}\}\,.
\]
This motivates the following:

\begin{conj}\label{conj5p4hu}
For every integer $n\ge 2$ and every field $k$ of characteristic different from $2$, one has
\[
p(F_{n-1}(x))=p(F_n)=\sup\{p(\ell(x_1,\dotsc, x_{n-1}))\,|\,\ell/k \text{ a finite field extension}\}\,.
\]
\end{conj}

By Prop.$\;$\ref{prop3p2hu}, the conjecture holds for arbitrary $n$ if $k$ is nonreal and everything is equal to $p(k(x_1,...,x_{n-1}))$ in this case. The conjecture also holds for general $k$ and $n=2$ due to Theorem$\;$\ref{thm3p1hu} (ii). In addition, we have shown that the conjecture holds for $n=3$ whenever $k$ is real and such that $p(k(x,y))\le 4$.

Moreover, if \cite[Conjecture$\;$4.15]{BGvG12}  were to hold for all real fields, then we could replace $\sup\{p(\ell(x_1,.... x_{n-1}))\}$ by $p(k(x_1,...,x_{n-1}))$ in the conjecture without making it weaker.

An immediate consequence of the conjecture is the inequality $p(\mathbb{R}(\!(t_1,\dotsc, t_{n})\!))\le 2^{n-1}$ for $n\ge 4$, which was conjectured in \cite{CDLR}.

\

For the $u$-invariant, similar considerations lead us to  propose  the following conjecture.

\begin{conj}\label{conj5p5hu}
For every integer $n\ge 2$ and every field $k$ of characteristic different from $2$, one has
\[
u(F_{n-1}(x))=u(F_n)=2\sup\{u(\ell(x_1,\dotsc, x_{n-1}))\,|\,\ell/k \text{ a finite field extension}\}\,.
\]
\end{conj}

For a real closed field $k$ and $n\ge 2$, if the equality $u(k(x_1,\dotsc, x_{n-1}))=2^{n-1}$ holds  as conjectured by Pfister \cite{Pfister82}, then the above conjecture implies $u(k(\!(t_1,\dotsc, t_n)\!))=2^{n}$.

\

\noindent {\bf Acknowledgements.} This work was done while the author was at the Universit\"at Duisburg-Essen and completed during his visit to the Max Planck Institute for Mathematics in Bonn. The author thanks Karim Becher for helpful discussions and comments. The author is grateful to the referee for the careful reading and many comments and suggestions, which significantly improved the exposition of the paper.

\addcontentsline{toc}{section}{\textbf{References}}

\bibliographystyle{alpha}
\bibliography{Laurent}

Author information:

\

Yong HU

\

Universit\'e de Caen, Campus 2

Laboratoire de Math\'ematiques Nicolas Oresme

14032, Caen Cedex

France

\end{document}